\documentclass[a4paper,11pt]{article}

\usepackage{amssymb,amsmath,amsthm}

\newtheorem{theorem}{Theorem}[section]
\newtheorem{lemma}[theorem]{Lemma}
\newtheorem{corollary}[theorem]{Corollary}
\newtheorem{conjecture}[theorem]{Conjecture}

\begin{document}
\title{Even pairs in square-free Berge graphs with no odd prism}

\author{Fr\'ed\'eric Maffray\thanks{CNRS, Laboratoire G-SCOP,
University of Grenoble-Alpes, France.  E-mail:
frederic.maffray@grenoble-inp.fr.}}

\date{\today}

 \maketitle

\begin{abstract}
We consider the class ${\cal G}$ of Berge graphs that contain no odd
prism and no square (cycle on four vertices).  We prove that every
graph $G$ in ${\cal G}$ either is a clique or has an even pair, as
conjectured by Everett and Reed.  This result is used to devise a
polynomial-time algorithm to color optimally every graph in ${\cal
G}$.
\end{abstract}

\noindent{\bf Keywords}: Berge graph, prism, square, even pair,
coloring, algorithm

\section{Introduction}

A graph $G$ is \emph{perfect} if every induced subgraph $H$ of $G$
satisfies $\chi(H)=\omega(H)$, where $\chi(H)$ is the chromatic number
of $H$ and $\omega(H)$ is the maximum clique size in $H$.  In a graph
$G$, a \emph{hole} is a chordless cycle with at least four vertices
and an \emph{antihole} is the complement of a hole.  Berge
{\cite{ber60,ber61,ber85}} introduced perfect graphs and conjectured
that a graph is perfect if and only if it does not contain as an
induced subgraph an odd hole or an odd antihole of length at least
$5$.  A \emph{Berge graph} is any graph that contains no odd hole and
no odd antihole of length at least $5$.  This famous question (the
Strong Perfect Graph Conjecture) was the objet of much research (see
\cite{ramree01}), until it was proved by Chudnovsky, Robertson,
Seymour and Thomas \cite{CRST}: \emph{Every Berge graph is perfect}.
Moreover, Chudnovsky, Cornu\'ejols, Liu, Seymour and Vu\v{s}kovi\'c
\cite{CCLSV} devised a polynomial-time algorithm that determines if a
graph is Berge (hence perfect).

It is known that one can obtain an optimal coloring of a perfect graph
in polynomial time due to the algorithm of Gr\"otschel, Lov\'asz and
Schrijver \cite{GLS}.  This algorithm however is not purely
combinatorial and impractical.  Here are some ideas that could be
fruitful in order to devise a purely combinatorial algorithm for
coloring perfect graphs.  An \emph{even pair} in a graph $G$ is a pair
$\{x,y\}$ of non-adjacent vertices such that every chordless path
between them has even length (number of edges).  Given two vertices
$x,y$ in a graph $G$, the operation of \emph{contracting} them means
removing $x$ and $y$ and adding one vertex with edges to every vertex
of $G\setminus \{x,y\}$ that is adjacent in $G$ to at least one of
$x,y$; we denote by $G/xy$ the graph that results from this operation.
Fonlupt and Uhry \cite{fonuhr82} proved that if $G$ is a perfect graph
and $\{x,y\}$ is an even pair in $G$, then the graph $G/xy$ is perfect
and $\chi(G/xy)=G$.  In particular, given a $\chi(G/xy)$-coloring $c$
of the vertices of $G/xy$, one can easily obtain a $\chi(G)$-coloring
of the vertices of $G$ as follows: keep the color for every vertex
different from $x,y$; assign to $x$ and $y$ the color assigned by $c$
to the contracted vertex.  This idea could be the basis for a
conceptually simple coloring algorithm for Berge graphs: as long as
the graph has an even pair, contract any such pair; when there is no
even pair find a coloring $c$ of the contracted graph and, applying
the procedure above repeatedly, derive from $c$ a coloring of the
original graph.

The algorithm for recognizing Berge graphs \cite{CCLSV} can be used to
detect an even pair in a Berge graph $G$; indeed, it is easy to see
that two non-adjacent vertices $a,b$ form an even pair in $G$ if and
only if the graph obtained by adding a vertex adjacent only to $a$ and
$b$ is Berge.  Thus, given a Berge graph $G$, one can try to color its
vertices by keeping contracting even pairs until none can be found.
Then some questions arise: which Berge graphs have no even pair, and
which do not?  What are the graphs for which a sequence of even-pair
contractions leads to graphs that are easy to color?

Bertschi \cite{ber90} proposed the following definitions.  A graph $G$
is \emph{even-contractile} if either $G$ is a clique or there exists a
sequence $G_0, \ldots, G_k$ of graphs such that $G=G_0$, for $i=0,
\ldots, k-1$ the graph $G_i$ has an even pair $\{x_i, y_i\}$ such that
$G_{i+1}=G_i/x_iy_i$, and $G_k$ is a clique.  A graph $G$ is
\emph{perfectly contractile} if every induced subgraph of $G$ is
even-contractile.  This class is of interest because it turns out that
many classical families of graphs are perfectly contractile; see
\cite{epsurvey}.  

Everett and Reed \cite{epsurvey} proposed a conjecture aiming at a
characterization of perfectly contractile graphs.  A \emph{prism} is a
graph that consists of two vertex-disjoint triangles (cliques of size
$3$) with three vertex-disjoint paths $P_1, P_2, P_3$ between them,
and with no other edge than those in the two triangles and in the
three paths.  The length of a path is its number of edges.  Note that
if two of the paths $P_1, P_2, P_3$ have lengths of different
parities, then their union induces an odd hole.  So in a Berge graph,
the three paths of a prism have the same parity.  A prism is
\emph{even} (resp.~\emph{odd}) if these three paths all have even
lengths (resp.~all have odd lengths).

\begin{conjecture}[\cite{epsurvey}]\label{conj1}
A graph is perfectly contractile if and only if it contains no odd 
hole, no antihole of length at least $5$, and no odd prism.
\end{conjecture}
Graphs that contain no odd hole, no antihole of length at least $5$,
and no odd prism were called \emph{Grenoble graphs} by Bruce Reed.

The `only if' part of Conjecture~ \ref{conj1} is not hard to
establish; see \cite{linmafree97} for the details.  The `if' part of
the conjecture remains open.  A weaker conjecture was proposed by
Everett and Reed \cite{epsurvey} and eventually proved by Maffray and
Trotignon \cite{MT}, as follows.

\begin{theorem}[\cite{MT}]\label{thm:MT}
If a graph contains no odd hole, no antihole of length at least $5$, 
and no prism then it is perfectly contractile.
\end{theorem}
The proof of Theorem~\ref{thm:MT} is a polynomial time algorithm that
takes as input any graph $G$ that contains no odd hole, no antihole of
length at least~$5$, and no prism, and produces a sequence of
contractions of even pairs that turns $G$ into a clique.  Moreover,
one can decide in polynomial time if a graph contains an odd hole, an
antihole of length at least~$5$ or a prism \cite{MTa}.

A \emph{square} is a hole of length four.  A graph is
\emph{square-free} if it does not contain a square as an induced
subgraph.  Here we will study Conjecture~\ref{conj1} in square-free
graphs.  We will be able to prove that every square-free Grenoble
graph that is not a clique has an even pair.  Unfortunately,
contracting an even pair may result in the presence of a square in the
contracted graph (if the two vertices of the even pair were linked by
a path of length four in the original graph).  So it is difficult to
establish that square-free Grenoble graphs are perfectly contractile.
Nevertheless, using the presence of even pairs, we will prove the
following theorem, which is the main result of this paper.

\begin{theorem}\label{thm:main}
There exists a combinatorial and polynomial time algorithm which,
given any square-free Grenoble graph $G$, returns an $\omega(G)$ coloring
of $G$ and a clique of size $\omega(G)$.
\end{theorem}

Since Theorem~\ref{thm:MT} settles the case of graphs that have no
prism, we may assume for our proof of Theorem~\ref{thm:main} that we
are dealing with a graph that contains an even prism.  So the next
sections focus on the study of such graphs.  Note that results
from~\cite{MTa} show that finding an induced prism in a Berge graph
can be done in polynomial time.

We finish this section with some notation and terminology.  In a graph
$G$, given a set $T\subset V(G)$, a vertex of $V(G) \setminus T$ is
\emph{complete to $T$} if it is adjacent to all vertices of $T$.  A
vertex of $V(G) \setminus T$ is \emph{anticomplete to $T$} if it is
not adjacent to any vertex of $T$.  Given two sets $S,T\subset V(G)$,
$S$ is \emph{complete} to $T$ if every vertex of $S$ is complete to
$T$, and $S$ is \emph{anticomplete} to $T$ if every vertex of $S$ is
anticomplete to $T$.  Given a path, any edge between two vertices that
are not consecutive along the path is a \emph{chord}.  A path that has
no chord is \emph{chordless}.

\section{Prisms}

Several sections in the proof of the Strong Perfect Graph
Theorem~\cite{CRST} are devoted to the analysis of Berge graphs that
contain a prism.  We extract here several theorems from~\cite{CRST}
that we will use.

Let $K$ be a prism, consisting of two vertex-disjoint triangles
$\{a_1, a_2, a_3\}$ and $\{b_1, b_2, b_3\}$, and three paths $P_1$,
$P_2$, $P_3$, where each $P_i$ has ends $a_i$ and $b_i$, and for $1\le
i < j \le 3$ the only edges between $V(P_i)$ and $V (P_j)$ are
$a_ia_j$ and $b_ib_j$.  The three paths $P_1$, $P_2$, $P_3$ are said
to \emph{form} the prism.  Vertices $a_1, a_2, a_3$ and $b_1, b_2,
b_3$ are the \emph{corners} of the prism.

\begin{theorem}[(7.3) in \cite{CRST}]\label{spgt73}
In a Berge graph $G$, let $R_1,R_2,R_3$ be three chordless paths that
form a prism $K$ with triangles $\{a_1, a_2, a_3\}$ and $\{b_1, b_2,
b_3\}$, where each $R_i$ has ends $a_i$ and $b_i$.  Assume that $R_1,
R_2, R_3$ all have length at least $2$.  Let $Y \subset V(G)$ be
anticonnected such that every vertex in $Y$ is adjacent to at least
two of $a_1, a_2, a_3$ and to at least two of $b_1, b_2, b_3$.  Then
at least two of $a_1, a_2, a_3$ and at least two of $b_1, b_2, b_3$
are complete to $Y$.
\end{theorem}

\begin{theorem}[(7.4) in \cite{CRST}]\label{spgt74}
In a Berge graph $G$, let $R_1,R_2,R_3$ be three chordless paths that
form a prism $K$ with triangles $\{a_1, a_2, a_3\}$ and $\{b_1, b_2,
b_3\}$, where each $R_i$ has ends $a_i$ and $b_i$.  Assume that $R_1,
R_2, R_3$ all have length at least $2$.  Let $R'_1$ be a chordless
path from $a'_1$ to $b_1$, such that $R'_1 , R_2, R_3$ also form a
prism.  Let $y \in V(G)$ have at least two neighbours in $A$ and in
$B$.  Then $y$ also has at least two neighbours in $\{a'_1, a_2,
a_3\}$.
\end{theorem}

\begin{theorem}[(10.1) in \cite{CRST}]\label{spgt101}
In a Berge graph $G$, let $R_1,R_2,R_3$ be three chordless paths that
form a prism $K$ with triangles $\{a_1, a_2, a_3\}$ and $\{b_1, b_2,
b_3\}$, where each $R_i$ has ends $a_i$ and $b_i$.  Let $F \subseteq
V(G) \setminus V(K)$ be connected, such that its set of attachments in
$K$ is not local.  Assume no vertex in $F$ is major with respect to
$K$.  Then there is a path $f_1$-$\ldots$-$f_n$ in $F$ with $n \ge 1$,
such that (up to symmetry) either:
\begin{itemize}
\item[1.]  
$f_1$ has two adjacent neighbours in $R_1$, and $f_n$ has two adjacent
neighbours in $R_2$, and there are no other edges between $\{f_1,
\ldots , f_n\}$ and $V(K)$, and (therefore) $G$ has an induced
subgraph which is the line graph of a bipartite subdivision of $K_4$,
or 
\item[2.]  
$n \ge 2$, $f_1$ is adjacent to $a_1, a_2, a_3$, and $f_n$ is adjacent
to $b_1, b_2, b_3$, and there are no other edges between $\{f_1,
\ldots, f_n\}$ and $V(K)$, or  
\item[3.]  
$n \ge 2$, $f_1$ is adjacent to $a_1, a_2$, and $f_n$ is adjacent to
$b_1, b_2$, and there are no other edges between $\{f_1, \ldots,
f_n\}$ and $V(K)$, or 
\item[4.] 
$f_1$ is adjacent to $a_1, a_2$, and there is at least one edge
between $f_n$ and $V(R_3) \setminus \{a_3\}$, and there are no other
edges between $\{f_1, \ldots, f_n\}$ and $V(K) \setminus \{a_3\}$.
\end{itemize}
\end{theorem}
In this paper the above theorem will always be applied to graphs that
do not contain any odd prism and (consequently) do not contain the
line-graph of any bipartite subdivision of $K_4$.  So only items~2, 3
or~4 hold.  Moreover, it is not specified that the prism is even in
the preceding theorem.  We will use the following special case of this
theorem.

\begin{corollary} 
In a Berge graph $G$, let $R_1,R_2,R_3$ be three chordless paths that
form a prism $K$ with triangles $\{a_1, a_2, a_3\}$ and $\{b_1, b_2,
b_3\}$, where each $R_i$ has ends $a_i$ and $b_i$ and has even length.
Let $x$ be a vertex in $V(G) \setminus V(K)$ such that $x$ is not a
major neighbor of $K$ and its set of attachments in $K$ is not local.
Then (up to symmetry) $x$ is adjacent to $a_1, a_2$, and there is at
least one edge between $x$ and $V(R_3) \setminus \{a_3, b_3\}$, and
there are no other edges between $x$ and $V(K) \setminus \{a_3\}$.
(In particular, $x$ is anticomplete to $\{b_1, b_2, b_3\}$.)
\end{corollary}

\begin{theorem}[(10.3) in \cite{CRST}]\label{spgt103}
Let $G$ be a Berge graph, such that there is no nondegenerate
appearance of $K_4$ in $G$.  Let $R_1,R_2,R_3$ form a prism $K$ in
$G$, with triangles $\{a_1, a_2, a_3\}$ and $\{b_1, b_2, b_3\}$, where
each $R_i$ has ends $a_i$ and $b_i$.  Let $F \subseteq V(G) \setminus
V(K)$ be connected, such that no vertex in $F$ is major with respect
to $K$.  Let $x_1$ be an attachment of $F$ in the interior of $R_1$,
and assume that there is another attachment $x_2$ of $F$ not in $R_1$.
Then there is a path $f_1$-$\ldots$-$f_n$ in $F$ such that (up to the
symmetry between $A$ and $B$) $f_1$ is adjacent to $a_2, a_3$, and
$f_n$ has at least one neighbour in $R_1 \setminus a_1$, and there are
no other edges between $\{f_1, \ldots, f_n\}$ and $V(K) \setminus
\{a_1\}$.
\end{theorem}

\section{Hyperprisms}

From now on, let $G$ be a square-free Berge graph that contains an
even prism.

We define hyperprisms as in \cite{CRST}.  Since $G$ contains an even
prism, $V(G)$ contains nine subsets
\begin{center}
\begin{tabular}{ccc}
$A_1$ & $C_1$ & $B_1$ \\ 
$A_2$ & $C_2$ & $B_2$ \\ 
$A_3$ & $C_3$ & $B_3$
\end{tabular}
\end{center}
with the following properties:
\begin{itemize}
\item
These nine sets are nonempty and pairwise disjoint.
\item
For distinct $i,j\in\{1,2,3\}$, $A_i$ is complete to $A_j$, and $B_i$
is complete to $B_j$, and there are no other edges between $A_i \cup
B_i \cup C_i$ and $A_j \cup B_j \cup C_j$.
\item
For each $i\in\{1,2,3\}$, every vertex of $A_i \cup B_i \cup C_i$
belongs to a chordless path between $A_i$ and $B_i$ with interior in
$C_i$.
\end{itemize}
The $9$-tuple $(A_1, C_1, B_1, A_2, C_2, B_2, A_3, C_3, B_3)$ is
called a \emph{hyperprism}.  For each $i\in\{1,2,3\}$, a chordless
path from $A_i$ to $B_i$ with interior in $C_i$ is called an
\emph{$i$-rung}.  Let us write $A = A_1 \cup A_2 \cup A_3$, $B = B_1
\cup B_2 \cup B_3$ and $C=C_1\cup C_2\cup C_3$.  Let $S_i = A_i \cup
B_i \cup C_i$ for $i\in\{1,2,3\}$.  The triple $(A_i, C_i, B_i)$ is
called a \emph{strip} of the hyperprism.  We call $(A,C,B)$ the
\emph{profile} of the hyperprism.

If we pick any $i$-rung $R_i$ for each $i\in\{1,2,3\}$, we see that
$R_1, R_2, R_3$ form a prism; any such prism is called an
\emph{instance} of the hyperprism.  Since $G$ contains no odd prism,
every instance of the hyperprism is an even prism, and so every rung
has even length.

Given two hyperprisms $\eta$ and $\eta'$ with profiles $(A,C,B)$ and
$(A',C',B')$ respectively, we write $\eta\prec \eta'$ if $C\subseteq
C'$ and either $A\subseteq A'$ and $B\subseteq B'$ or $A\subseteq B'$
and $B\subseteq A'$ and one of these inclusions is strict.  Clearly,
$\prec$ is an order relation on hyperprisms, so we can speak about
maximal hyperprisms for $\prec$.  Although the notion of profile does
not appear in \cite{CRST}, it is easy to see that the notion of
maximal hyperprism in \cite{CRST} is equivalent to that which is
defined here.

Let $\eta=(A_1, \ldots, B_3)$ be a hyperprism, and let $H$ be the
subgraph of $G$ induced on the union of these nine sets.  A subset $X
\subseteq V(H)$ is \emph{local} (with respect to the hyperprism) if
$X$ is a subset of one of $S_1$, $S_2$, $S_3$, $A$ or $B$.  Let $x$ be
any vertex in $V(G)\setminus V(H)$.  We say that $x$ is a \emph{major}
neighbor of $H$ is $x$ is a major neighbor of an instance of $H$.  Let
$M$ be the set of all major neighbors of $H$.

From now on, we assume that $\eta$ is a maximal hyperprism.

\begin{lemma}\label{lem:spgt2106}
For every connected subset $F$ of $V(G) \setminus (V(H)\cup M)$, its
set of attachments in $H$ is local.
\end{lemma}
This lemma is identical to Claim~(2) in the proof of Theorem~(10.6) in
\cite{CRST}, so we omit its proof.

\begin{lemma}\label{lem:mabcuts}
For each $i\in\{1,2,3\}$, $M\cup A_i\cup B_i$ is a cutset that
separates $C_i$ from $S_{i+1}\cup S_{i+2}$.  Consequently, $C_1$,
$C_2$ and $C_3$ lie in three distinct components of $G\setminus (M\cup
A\cup B)$.
\end{lemma}
\noindent{\it Proof.} For suppose on the contrary that there is a path
$P=p$-$\cdots$-$q$, with $V(P)\subset V(G)\setminus (M\cup A_i\cup
B_i)$ such that $p$ has a neighbor in $C_i$ and $q$ has a neighbor in
$S_{i+1}\cup S_{i+2}$.  Let $P$ be a shortest such path; then
$V(P)\subseteq V(G)\setminus V(H)$, so $P$ contradicts
Lemma~\ref{lem:spgt2106}.  \hfill $\Box$

\begin{lemma}\label{lem:xaabb}
Let $x\in M$.  Then $x$ is complete to at least two of $A_1, A_2, A_3$
and at least two of $B_1, B_2, B_3$.
\end{lemma}
\noindent{\it Proof.} Since $x$ is in $M$, there exists for each
$i\in\{1, 2, 3\}$ an $i$-rung $R_i$ such that $x$ is a major neighbor
of the prism $K$ formed by $R_1, R_2, R_3$.  Let $R_i$ have ends
$a_i\in A_i$ and $b_i\in B_i$ ($i=1,2,3$).  Consider any $1$-rung
$P_1$, and let $K'$ be the prism formed by $P_1, R_2, R_3$.  We claim
that:
\begin{equation}
\mbox{$x$ is a major neighbor of $K'$.}\label{eq:xmkp}
\end{equation}
For suppose the contrary.  Let $X$ be the set of neighbors of $x$.
Let $P_1$ have ends $a'_1\in A_1$ and $b'_1\in B_1$, and let
$A'=\{a'_1, a_2, a_3\}$ and $B'=\{b'_1, b_2, b_3\}$.  If $b'_1=b_1$,
then Theorem~\ref{spgt74} shows that $x$ has at least two neighbors in
$A'$, and so the claim holds.  Therefore assume that $b'_1\neq b_1$
and, similarly, that $a'_1\neq a_1$.  Let $\alpha=|X\cap A|$,
$\beta=|X\cap B|$, $\alpha'=|X\cap A'|$, $\beta'=|X\cap B'|$.  We know
that $\alpha\ge 2$ and $\beta\ge 2$ since $x$ is a major neighbor of
$K$, and $\min\{\alpha',\beta'\}\le 1$ since $x$ is not a major
neighbor of $K'$.  Moreover, $\alpha'\ge \alpha-1$ and $\beta'\ge
\beta-1$ since $K$ and $K'$ differ by only one rung.  Up to the
symmetry on $A,B$, these conditions imply that the vector
$(\alpha,\beta,\alpha',\beta')$ is equal to either $(3,2,3,1)$,
$(3,2,2,1)$, $(2,2,2,1)$ or $(2,2,1,1)$.  In either case we have
$\beta=2$ and $\beta'=1$, so $x$ is adjacent to $b_1$, not adjacent to
$b'_1$, and adjacent to exactly one of $b_2, b_3$, say to $b_3$.  \\
Suppose that $(\alpha',\beta')$ is equal to $(3,1)$ or $(2,1)$.  We
can apply Theorem~\ref{spgt101} to $K'$ and $F=\{x\}$, and it follows
that $x$ satisfies item~4 of that theorem, so $x$ is adjacent to
$a'_1, a_2, b_3$ and has no neighbor in $V(K') \setminus (\{a'_1,
a_2\}\cup V(R_3))$.  In particular $x$ has no neighbor in
$V(R_2)\setminus\{a_2\}$, and then $V(R_2)\cup\{x,b_3\}$ induces an
odd hole, a contradiction.  So we may assume that
$(\alpha,\beta,\alpha',\beta') = (2,2,1,1)$, which restores the
symmetry between $A$ and $B$.  Since $\alpha=2$ and $\alpha'=1$, $x$
is adjacent to $a_1$, not adjacent to $a'_1$, and adjacent to exactly
one of $a_2, a_3$.  In fact if $x$ is adjacent to $a_2$, then $K'$ and
$\{x\}$ violate Theorem~\ref{spgt101}.  So $x$ is adjacent to $a_3$
and not to $a_2$, and Theorem~\ref{spgt101} implies that $x$ is a
local neighbor of $K'$ with $X\cap K'\subseteq V(R_3)$, so $x$ has no
neighbor on $P_1$ or $R_2$.  \\
We observe that for every $1$-rung $Q_1$, the ends of $Q_1$ are either
both adjacent to $x$ or both not adjacent to $x$, for otherwise the
prism formed by $Q_1, R_2, R_3$ and the set $F=\{x\}$ violate
Theorem~\ref{spgt101}.  Let $A'_1=A_1\setminus X$ and $A''_1=A_1\cap
X$, and similarly $B'_1=B_1\setminus X$ and $B''_1=B_1\cap X$.  The
preceding observation means that every $1$-rung is either between
$A'_1$ and $B'_1$ or between $A''_1$ and $B''_1$.  Let $C'_1$ be the
set of vertices of $C_1$ that lie on a $1$-rung whose ends are in
$A'_1\cup B'_1$, and let $C''_1$ be the set of vertices of $C_1$ that
lie on a $1$-rung whose ends are in $A''_1\cup B''_1$.  The sets
$C'_1$ and $C''_1$ are disjoint and there is no edge between $A'_1\cup
C'_1\cup B'_1$ and $C''_1$ or between $A''_1\cup C''_1\cup B''_1$ and
$C'_1$, for otherwise we would find a $1$-rung with one end in $A'_1$
and the other in $B''_1$.  For every $1$-rung $P'_1$ with ends in
$A'_1\cup B'_1$ Theorem~\ref{spgt101} implies (just like for $P_1$)
that $x$ is a local neighbor of the prism formed by $P'_1, R_2, R_3$,
so $x$ has no neighbor on $P'_1$.  Hence $x$ has no neighbor in
$A'_1\cup C'_1\cup B'_1$.  We claim that $A'_1$ is complete to
$A''_1$.  For suppose on the contrary, up to relabelling vertices and
rungs, that $a'_1$ and $a_1$ are not adjacent.  Then $V(R_1)\cup\{x,
a_1, a_2, b_3\}$ induces an odd hole.  So the claim holds, and
similarly $B'_1$ is complete to $B''_1$.  \\
Now we consider $S_2$.  Let $A'_2=A_2\setminus X$, $A''_2=A_2\cap
X$, $B'_2=B_2\setminus X$ and $B''_2=B_2\cap X$.  By the same
arguments as for the $1$-rungs, we see that every $2$-rung $Q_2$ is
either between $A'_2$ and $B'_2$ or between $A''_2$ and $B''_2$, for
otherwise the prism formed by $P_1, Q_2, R_3$ and the set $F=\{x\}$
violate Theorem~\ref{spgt101}.  Let $C'_2$ be the set of vertices of
$C_2$ that lie on a $2$-rung whose ends are in $A'_2\cup B'_2$, and
let $C''_2$ be the set of vertices of $C_2$ that lie on a $1$-rung
whose ends are in $A''_2\cup B''_2$.  Then, by the same arguments as
above, $C'_2$ and $C''_2$ are disjoint and there is no edge between
$A'_2\cup C'_2\cup B'_2$ and $C''_2$ or between $A''_2\cup C''_2\cup
B''_2$ and $C'_2$.  Also $x$ has no neighbor in $A'_2\cup C'_2\cup
B'_2$, and $A'_2$ is complete to $A''_2$ and $B'_2$ is complete to
$B''_2$.  It follows that the nine sets
\begin{center}
\begin{tabular}{ccc}
$A'_1$ & $C'_1$ & $B'_1$ \\ 
$A'_2$ & $C'_2$ & $B'_2$ \\ 
$A''_1\cup A''_2\cup A_3$ & $C''_1\cup C''_2\cup C_3\cup\{x\}$ &
$B''_1\cup B''_2\cup B_3$
\end{tabular}
\end{center}
form a hyperprism, which contradicts the maximality of $\eta$.  Thus
(\ref{eq:xmkp}) holds.

By (\ref{eq:xmkp}) applied repeatedly, we obtain that $x$ is a major
neighbor of every instance of $H$.

Now suppose that $x$ has a non-neighbor $u_1\in A_1$ and a
non-neighbor $u_2\in A_2$.  For each $i\in\{1,2\}$ let $P_i$ be an
$i$-rung with end $u_i$, and let $P_3$ be any $3$-rung.  Then $x$ is
not a major neighbor of the prism formed by $P_1, P_2, P_3$, a
contradiction.  So $x$ is complete to one of $A_1, A_2$, say to $A_1$.
Likewise, $x$ is complete to one of $A_2, A_3$.  So $x$ is complete to
at least two of $A_1, A_2, A_3$.  The same holds for $B_1, B_2, B_3$.
This completes the proof of the lemma.  \hfill $\Box$

\begin{lemma}\label{lem:maabb}
Let $M$ be the set of major neighbors of~$\eta$.  Then: \\
(i) Two of $A_1, A_2, A_3$ and two of $B_1, B_2, B_3$ are cliques.  \\
(ii) $M$ is complete to at least two of $A_1, A_2, A_3$ and at least
two of $B_1, B_2, B_3$.  \\
(iii) There is an integer $j\in\{1, 2, 3\}$ such that $A_j$ and $B_j$
are cliques and $M$ is complete to $A_j\cup B_j$.
\end{lemma}
\noindent{\it Proof.} If (i) does not hold, then, up to symmetry,
there are two non-adjacent vertices in $A_1$ and two non-adjacent
vertices in $A_2$, and these four vertices induce a square, a
contradiction.

(ii) We claim that $M$ is complete to one of $A_1, A_2$.  For suppose
on the contrary that there are two non-adjacent vertices $a_1\in A_1$
and $u\in M$ and also two non-adjacent vertices $a_2\in A_2$ and $v\in
M$.  By Lemma~\ref{lem:xaabb}, $u$ is complete to $A_2$ and $v$ is
complete to $A_1$, so $ua_2$ and $va_1$ are edges, and $u\neq v$.  If
$u$ and $v$ are not adjacent, then, by Theorem~\ref{spgt73} applied to
$K$ and $Y=\{u,v\}$, there is a vertex $b\in B$ that is complete to
$Y$, and then $\{a_1, a_2, u, v, b\}$ induces a $5$-hole, a
contradiction.  So $u$ and $v$ are adjacent, and $\{u,v,a_1,a_2\}$
induces a square, a contradiction.  So the claim holds, say $M$ is
complete to $A_1$.  Similarly, $M$ is complete to one of $A_2, A_3$.
Thus $M$ is complete to two of $A_1, A_2, A_3$, and the same holds for
$B_1, B_2, B_3$ by symmetry.

(iii) By~(ii), we may assume that $M$ is complete to $A_1\cup B_1$.
If both $A_1, B_1$ are cliques, then (iii) holds with $j=1$.
Therefore assume that $A_1$ is not a clique.  By~(i), $A_2$ and $A_3$
are cliques.  Moreover $M$ is complete to $A_2\cup A_3$, for if there
are non-adjacent vertices $u\in M$ and $a\in A_2\cup A_3$, then by
Lemma~\ref{lem:xaabb} the vertex $u$ is complete to $A_1$, and then
$u,a$ and two non-adjacent vertices from $A_1$ induce a square.
By~(ii) $M$ is complete to one of $B_2, B_3$, say to $B_2$.  So if
$B_2$ is a clique, then (iii) holds with $j=2$.  Therefore assume that
$B_2$ is not a clique.  Then $B_3$ is a clique by~(i), moreover, as
above (with $A_1$), $M$ is complete to $B_3$.  So (iii) holds with
$j=3$.  Thus the lemma is proved.  \hfill $\Box$

\subsection{Selecting a strip}

Let us say that a strip $(A_i, C_i, B_i)$ of the hyperprism is
\emph{good} if both $A_i$ and $B_i$ are cliques and $M$ is complete to
$A_i\cup B_i$.  Lemma~\ref{lem:maabb} says that every maximal
hyperprism has a good strip.  We may assume that $(A_1, C_1, B_1)$ is
a good strip of $\eta$.  Moreover, we may assume that we choose $\eta$
such that $S_1$ has the smallest size over all good strips of maximal
hyperprisms.

\begin{lemma}\label{lem:path}
Let $P=a$-$u$-$\cdots$-$v$-$b$ be any chordless path with $a\in A_1$,
$b\in B_1$, and $V(P)\cap M=\emptyset$.  Then $V(P)\subset V(H)$.
Moreover, either: \\
$\bullet$ $P$ is a $1$-rung, or \\
$\bullet$ the interior of $P$ is an $i$-rung for some $i\in\{1,2,3\}$, or \\
$\bullet$ $P$ has odd length, $V(P)\subseteq S_1$ and exactly one of $u\in
A_1$ and $v\in B_1$ holds.
\end{lemma}
\noindent{\it Proof.} Note that $P$ has length at least~$2$.  We prove
the lemma by induction on the length of $P$.  If $P$ has length~$2$,
say $P=a$-$x$-$b$, then we must have $x\in C_1$ (for otherwise, we
could add $x$ to $C_1$ and obtain a hyperprism that contradicts the
maximality of $\eta$), and so $P$ is a $1$-rung.  Now assume that the
length of $P$ is at least~$3$.  Let $\tilde{P}$ be the interior of
$P$.

When $V(P)\not\subset V(H)$, there are subpaths $P_1, \ldots, P_k$ of
$P$, with $k$ odd, $k\ge 3$, such that $P=P_1$-$P_2$-$\cdots$-$P_k$,
with $a\in V(P_1)$ and $b\in V(P_k)$, and, for all odd $j$,
$V(P_j)\subset V(H)$, and for all even $j$, $V(P_j)\cap
V(H)=\emptyset$.  When $V(P)\subset V(H)$ we use the same notation,
with $k=1$.  When $k\ge 3$, for each even $j$, let $X_j$ be the set of
attachment of $P_j$ in $H$.  We claim that:
\begin{equation}\label{xjplocal}
\mbox{For each even $j$, there is $i_j\in\{1,2,3\}$ such that
$X_j\subseteq S_{i_j}$.}
\end{equation}
By Lemma~\ref{lem:spgt2106} applied to $P_j$, we know that $X_j$ is
local with respect to $H$.  Suppose that $X_j\subseteq A$.  Let $w$
(resp.~ $w'$) be the vertex in $P_{j-1}$ (resp.~in $P_{j+1}$) that has
a neighbor in $P_j$.  Then $w,w'\in X_j$, and $w,w'$ are not adjacent,
so $w,w'\in A_2\cup A_3$ and $a$ is adjacent to both $w,w'$, a
contradiction.  Hence $X_j$ is not a subset of $A$ and, similarly, not
of $B$ either.  Thus (\ref{xjplocal}) holds.

\medskip

Suppose that $u\in A_2$.  Then $V(P_1)\cap A=\{a,u\}$ (for otherwise
$a$ would have a neighbor on $P\setminus \{a, u\}$), and so
$V(P_1\setminus a)\subseteq S_2$.  Now, applying (\ref{xjplocal})
repeatedly, we obtain that for each even $j$ we have $X_j\subseteq
S_2$, for each odd $j$ with $j<k$ we have $V(P_j)\subseteq S_2$, and
$V(P_k \setminus b)\subseteq S_2$.  Then $V(\tilde{P})\subseteq S_2$,
for otherwise we could add the vertices of $\tilde{P}$ to $S_2$ and
thus obtain a hyperprism that contradicts the maximality of $\eta$.
Hence $\tilde{P}$ is a $2$-rung and the lemma holds.  We obtain a
similar conclusion if either $u\in A_3$ or $v\in B_2\cup B_3$.  Now
assume that $u\notin A_2\cup A_3$ and $v\notin B_2\cup B_3$.

Suppose that $u\notin A_1$ and $v\notin B_1$.  Then $V(P_1)\cap
A=\{a\}$ and $V(P_1)\subseteq S_1$.  Now, applying (\ref{xjplocal})
repeatedly, we obtain that for each even $j$ we have $X_j\subseteq
S_1$, for each odd $j$ with $j<k$ we have $V(P_j)\subseteq S_1$, and
$V(P_k\setminus b)\subseteq S_1$.  Then $V(P)\subseteq S_1$, for
otherwise we could add the vertices of $P$ to $S_1$ and thus obtain a
hyperprism that contradicts the maximality of $\eta$.  Hence $P$ is a
$1$-rung and the lemma holds.

Now suppose that $u\in A_1$ and $v\in B_1$.  We can apply induction to
$\tilde{P}$.  It cannot be that the second or third item of the lemma
holds for $\tilde{P}$ (for otherwise $a$ would have two neighbors on
$P$), so the first item holds for $\tilde{P}$, and so the second item
holds for $P$.

Finally suppose, up to symmetry, that $u\in A_1$ and $v\notin B_1$.
We can apply induction to $P\setminus a$.  It cannot be that the
second or third item of the lemma holds for $P\setminus a$, so the
first item holds for $P\setminus a$, and so the third item holds for
$P$.  Thus the lemma holds.  \hfill $\Box$

\medskip

A \emph{necklace} is a graph that consists of four disjoint chordless
paths $R_1= a\cdots a'$, $R_2=b\cdots b'$, $R_3=c\cdots c'$, $R_4=
d\cdots d'$, where $R_1, R_2$ may have length $0$ but $R_3$, $R_4$
have length at least $1$, and such that the edge-set of $S$ is
$E(R_1)\cup E(R_2)\cup E(R_3) \cup E(R_4)\cup\{a'c, a'd, cd, b'c',
b'd', c'd'\}$.  Note that $\{a',c,d\}$ and $\{b',c',d'\}$ are
triangles in $S$.  Vertices $a$ and $b$ are the endvertices of the
necklace, and we may also say that $S$ is an $(a, b)$-necklace.

\medskip

Let $R'$ and $R''$ be two $1$-rungs, where $R'$ has ends $u',w$, and
$R''$ has ends $u'',w$, and $u'\neq u''$ (so $w$ is in one of the two
sets $A_1, B_1$ and $u',u''$ are in the other set).  We say that $R'$
and $R''$ \emph{converge} if $u'$ has no neighbor in $R''\setminus
u''$ and $u''$ has no neighbor in $R'\setminus u'$.

\begin{lemma}\label{lem:conv}
There do not exist two $1$-rungs that converge.
\end{lemma}
\noindent{\it Proof.} Suppose on the contrary that $R'$ and $R''$ are
two $1$-rungs that converge.  Choose $R'$ and $R''$ such that
$|V(R')\cup V(R'')|$ is minimized.  Let $R'=u_0$-$u_1$-$\cdots$-$u_p$
(with $p$ even, $p\ge 2$) and $R''=v_0$-$v_1$-$\cdots$-$v_q$ (with $q$
even, $q\ge 2$), and assume up to symmetry that $u_0, v_0\in A_1$,
$u_0\neq v_0$, $u_p=v_q\in B_1$, $u_0$ has no neighbor in
$R''\setminus v_0$, and $v_0$ has no neighbor in $R'\setminus u_0$.
Let $i$ be the smallest integer such that $u_i$ has a neighbor in
$R''\setminus v_0$.  Note that $i$ exists since $u_{p-1}$ has a
neighbor in $R''\setminus v_0$.  Also $i\neq 0$ because of the
hypothesis on $u_0$.  Likewise, let $j$ be the smallest integer such
that $v_j$ has a neighbor in $R'\setminus u_0$.  Let $h$ be the
smallest integer such that $u_iv_h$ is an edge.  So $0<j\le h$.
Moreover, $h<q$, for otherwise we must have $i=p-1$ and $V(R')\cup
V(R'')$ induces an odd hole.  Now the set $\{u_0, \ldots, u_i, v_0,
\ldots, v_h\}$ induces a hole, so it is an even hole, so $i$ and $h$
have the same parity.  We claim that:
\begin{equation}\label{ipjp}
\mbox{We may assume that $R'[u_{i+1}, u_p]= R''[v_{j+1}, v_q]$.}
\end{equation}
To prove this, first suppose that $i\neq p-1$.  Let $k$ be the largest
integer such that $u_iv_k$ is an edge.  Then
$u_0$-$u_1$-$\cdots$-$u_i$-$v_k$-$\cdots$-$v_q$ is a chordless path,
so it is a $1$-rung, and it must have even length, so $h$ and $k$ have
different parities.  If $k\neq h+1$, then
$v_0$-$v_1$-$\cdots$-$v_h$-$u_i$-$v_k$-$\cdots$-$v_q$ is a chordless
path, so it is a $1$-rung, and it has odd length, a contradiction.
Hence $k=h+1$.  The minimality of $|V(R')\cup V(R'')|$ implies that
$R'[u_{i+1}, u_p]= R''[v_{h+1}, v_q]$, so $h=j$ and the claim holds.
Therefore we may assume that $i=p-1$.  By the same argument as with
$i$, we may assume that $j=q-1$ (so $h=q-1$).  Thus (\ref{ipjp})
holds.

Let $R_2$ be any $2$-rung, with ends $a_2\in A_2$ and $b_2\in B_2$.
Let $P_1=u_0$-$u_1$-$\cdots$-$u_i$, $P_2=v_0$-$v_1$-$\cdots$-$v_j$,
and $P_3=a_2$-$R_2$-$b_2$-$u_p$-$u_{p-1}$-$\cdots$-$u_{i+1}$.  It
follows from (\ref{ipjp}) that $P_1, P_2, P_3$ form a prism.  Since
$G$ contains no odd prism, these three paths have even length, and so
$i$ and $j$ are even, so $u_{i+1}\neq u_p$ and $u_{i+1}\in C_1$.  Let
$Z=C_2\cup C_3\cup B_2\cup B_3\cup\{u_{i+2}, \ldots, u_p\}$.  We
observe that:
\begin{center}
\begin{tabular}{ccc}
$\{u_0\}$ & $\{u_1, \ldots, u_{i-1}\}$ & $\{u_i\}$ \\ 
$\{v_0\}$ & $\{v_1, \ldots, v_{i-1}\}$ & $\{v_j\}$ \\ 
$A_2\cup A_3$ & $Z$ & $\{u_{i+1}\}$
\end{tabular}
\end{center}
form a hyperprism $\eta'$.  So there exists a maximal hyperprism
$\eta^*$ such that $\eta'\preceq \eta^*$.  Let $\eta^*=(A_1^*, C_1^*,
B_1^*, A_2^*, C_2^*, B_2^*, A_3^*, C_3^*, B_3^*)$, $A^*=A_1^*\cup
A_2^*\cup A_3^*$, $B^*=B_1^*\cup B_2^*\cup B_3^*$ and $C^*=C_1^*\cup
C_2^*\cup C_3^*$, and, for each $i\in\{1,2,3\}$, $S_i^*=A_i^*\cup
C_i^*\cup B_i^*$.  We know that $\{u_0, v_0\}\cup A_2\cup A_3\subseteq
A^*$, and $\{u_i, v_j, u_{i+1}\}\subseteq B^*$, and $\{u_1, \ldots,
u_{i-1}\}\cup \{v_1, \ldots, v_{i-1}\}\cup Z\subseteq C^*$.  Since $Z$
is connected, we may assume, up to symmetry, that $Z\subseteq C_3^*$,
and so $A_2\cup A_3\subseteq A_3^*$ and $\{u_{i+1}\}\subseteq B_3^*$.
We claim that:
\begin{equation}\label{S12starinS1}
\mbox{$S_1^*\cup S_2^*\subset S_1$, and $A_1^*\cup A_2^*\subset A_1$
and $B_1^*\cup B_2^*\cup B_3^*\subset C_1$.}
\end{equation}
We may assume up to symmetry that $P_2$ is either a $2$-rung or a
$3$-rung of $\eta^*$.  Let $R_1^*$ be any $1$-rung of $\eta^*$, with
ends $a_1^*\in A_1^*$ and $b_1^*\in B_1^*$.  So $a_1^*$ is complete to
$A_2\cup A_3\cup\{v_0\}$ and $b_1^*$ is complete to $\{v_j,
u_{i+1}\}$, and there are no other edges between $V(R_1^*)$ and
$V(P_2)\cup V(P_3)$.  Let
$R_1=a_1^*$-$R_1^*$-$b_1^*$-$u_{i+1}$-$R'$-$u_p$; so $R_1$ is an even
chordless path.  Let $R_1^+=v_0$-$a_1^*$-$R_1$-$u_p$; so $R_1^+$ is an
odd chordless path.  By Lemma~\ref{lem:path}, we have $a_1^*\in A_1$
and $R_1$ is a $1$-rung of $\eta$.  Thus $V(R_1^*)\subset A_1\cup C_1$
for every $1$-rung $R_1^*$ of $\eta^*$, and $A_1^*\subset A_1$ and
$B_1^*\subset C_1$.  We see that $R_1$ converges with $R''$, so we may
let $R_1^*$ play the role of $R'$, which restores the symmetry between
$1$-rungs and $2$-rungs of $\eta^*$, and consequently $V(R_2^*)\subset
S_2$ holds for every $2$-rung $R_2^*$ of $\eta^*$.  Thus
(\ref{S12starinS1}) holds.

Let $M^*$ be the set of major neighbors of $\eta^*$.  We claim that:
\begin{equation}\label{eq:mstara12}
\mbox{$M^*$ is complete to $A_1^*\cup A_2^*$.}
\end{equation}
For suppose that some vertex $m^*\in M^*$ is not is complete to
$A_1^*\cup A_2^*$.  Then $m^*$ is complete to $A_3^*$ and in
particular to $A_2\cup A_3$.  Moreover $m^*\notin A_1$, since $A_1$ is
a clique and $A_1^*\cup A_2^*\subset A_1$.  Therefore $m^*\notin
V(H)$.  We know that $m^*$ is complete to one of $B_1^*, B_2^*$, which
are subsets of $C_1$.  Hence the set of attachments of $m^*$ in $H$ is
not local, so Lemma~\ref{lem:spgt2106} implies that $m^*\in M$, so
$m^*$ is complete to $A_1$, a contradiction.  Thus (\ref{eq:mstara12})
holds.

\begin{equation}\label{Sjstargood}
\mbox{For some $j\in\{1,2\}$, $(A_j^*, C_j^*, B_j^*)$ is a good strip 
of $\eta^*$.}
\end{equation}
Since $A_1^*\cup A_2^*\subseteq A_1$, both $A_1^*$ and $A_2^*$ are
cliques.  We may assume up to symmetry that $M^*$ is complete to
$B_1^*$.  So if $B_1^*$ is a clique, the claim holds with $j=1$.  Now
assume that $B_1^*$ is not a clique.  Then $B_2^*$ is a clique by
Lemma~\ref{lem:maabb} applied to $\eta^*$, and $M^*$ is complete to
$B_2^*$ (for otherwise two non-adjacent vertices from $M^*\cup B_2^*$
plus two non-adjacent vertices from $B_1^*$ induce a square, and so
the claim holds with $j=2$.  Thus (\ref{Sjstargood}) holds.
 
Now Claims~(\ref{S12starinS1}) and (\ref{Sjstargood}) contradict the
choice of $\eta$ (with the smallest good strip).  This completes the
proof of the lemma.  \hfill $\Box$

\subsection{Finding an even pair}

Pick any $b\in B_1$.  For any two $a,a'\in A_1$, write $a <_b a'$
whenever there exists an odd chordless path $R$ from $a$ to $b$ such
that $a'$ is the neighbor of $a$ on $R$.  Note that in that case,
Lemma~\ref{lem:path} implies that $R\setminus a$ is a $1$-rung.

\begin{lemma}\label{lem:order}
For each $b\in B_1$, $<_b$ is an order relation.
\end{lemma}
\noindent{\it Proof.} We first claim that the relation $<_b$ is
antisymmetric.  Suppose on the contrary that there are vertices
$u,v\in A_1$ such that $u<_b v$ and $v<_b u$.  So there exists an odd
chordless path $P_u=u$-$v$-$\cdots$-$b$ and there exists an odd
chordless path $P_v=v$-$u$-$\cdots$-$b$.  By Lemma~\ref{lem:path},
$P_u\setminus u$ and $P_v\setminus v$ are $1$-rungs.  Because of $b$
these two rungs converge, which contradicts Lemma~\ref{lem:conv}.  So
$<_b$ is antisymmetric.  
Now we claim that $<_b$ is transitive.  Let $u,v,w$ be three vertices
in $A_1$ such that $u<_b v<_b w$.  So there is an odd chordless path
$v$-$w_0$-$w_1$-$\cdots$-$w_k$ with $k$ even, $k\ge 2$, $w=w_0$ and
$w_k=b$.  By Lemma~\ref{lem:path}, $w_0$-$w_1$-$\cdots$-$w_k$ is a
$1$-rung.  Let $j$ be the largest integer such that $uw_j$ is an edge.
Suppose that $j>0$.  If $j$ is even, then $u$-$w_j$-$\cdots$-$w_k$ is
a $1$-rung of odd length, a contradiction.  If $j$ is odd, then
$v$-$u$-$w_j$-$\cdots$-$w_k$ is an odd chordless path, so $v<_b u$,
which contradicts the fact that $<_b$ is antisymmetric.  So $j=0$,
which implies that $u<_b w$.  Hence $<_b$ is antisymmetric and
transitive, so it is an order relation.  \hfill $\Box$

\medskip

Similarly, for each $a\in A_1$, and for any two $b,b'\in B_1$, we
write $b <_a b'$ whenever there exists an odd chordless path $R$ from
$b$ to $a$ such that $b'$ is the neighbor of $b$ on $R$.  So $<_a$ is
an order relation on $B_1$ for each $a$.  

\begin{lemma}\label{lem:twist}
If there are four vertices $a,u\in A_1$ and $b,v\in B_1$ such that $a 
<_b u$ and $b <_u v$, then $a <_v u$.
\end{lemma}
\noindent{\it Proof.} The hypothesis that $a<_b u$ means that there is
an odd chordless path $R=a$-$r_0$-$r_1$-$\cdots$-$r_k$ with $r_0=u$
and $r_k=b$.  By Lemma~\ref{lem:path}, $R\setminus a$ is a $1$-rung,
so $k$ is even.  The hypothesis that $b <_u v$ means that there is an
odd chordless path $Q=b$-$v$-$\cdots$-$u$, and, by
Lemma~\ref{lem:path}, $Q\setminus b$ is a $1$-rung.  If $v$ has no
neighbor in $R\setminus b$, then $R\setminus a$ and $Q\setminus b$ are
two rungs that converge, a contradiction.  So there is an integer
$j<k$ such that $vr_j$ is an edge, and we choose the smallest such
$j$.  So $r_0$-$r_1$-$\cdots$-$r_j$-$v$ is a $1$-rung, so $j$ is odd.
Then $a$-$r_0$-$r_1$-$\cdots$-$r_j$-$v$ is an odd chorldess path,
which shows that $a <_v r_0$, i.e., $a <_v u$.  \hfill $\Box$

\begin{lemma}\label{lem:a1b1}
There exists an even pair $\{a,b\}$ with $a\in A_1$ and $b\in B_1$.
\end{lemma}
\noindent{\it Proof.} For each $a\in A_1$, let $\mbox{Max}(a)$ be the
set of maximal elements of the partially ordered set $(B_1, <_a)$.
Likewise, for each $b\in B_1$, let $\mbox{Max}(b)$ be the set of
maximal elements of $(A_1, <_b)$.  We claim that:
\begin{equation}\label{maxab}
\mbox{There exist $a\in A_1$ and $b\in B_1$ such that
$a\in\mbox{Max}(b)$ and $b\in\mbox{Max}(a)$.}
\end{equation}
For each $a\in A_1$ and $b\in B_1$, let $D(a,b)=\{a'\in A \mid a' <_b
a\}$.  Choose $a$ and $b$ such that the size of $D(a,b)$ is maximized.
We have $a\in\mbox{Max}(b)$, for otherwise, there is $u\in A_1$ such
that $a <_b u$, so $D(u,b)\supseteq D(a,b)\cup\{a\}$, which
contradicts the choice of $a$ and $b$.  So if $b\in\mbox{Max}(a)$ the
claim holds.  Hence let us assume that $b\notin\mbox{Max}(a)$.  This
means that there exists $v\in\mbox{Max}(a)$ such that $b <_a v$.  If
$a\in\mbox{Max}(v)$, then the claim holds with the pair $a,v$.  Hence
let us assume that $a\notin\mbox{Max}(v)$.  So there exists
$u\in\mbox{Max}(v)$ such that $a <_v u$.  For each $a'\in D(a,b)$, we
can apply Lemma~\ref{lem:twist} to the four vertices $a',a,b,v$, which
implies $a' <_v a$ and (by the transitivity of $<_v$) $a' <_v u$.  So
$D(u,v)\supseteq D(a,b)\cup\{a\}$, which contradicts the choice of $a$
and $b$.  Thus (\ref{maxab}) holds.

Let $a,b$ be any two vertices that satisfy (\ref{maxab}).  We claim
that $\{a,b\}$ is an even pair of $G$.  For suppose that there exists
an odd chordless path $P$ with ends $a$ and $b$.  By
Lemma~\ref{lem:path}, and up to symmetry, we may assume that the
neighbor $a'$ of $a$ on $P$ is in $A_1$ and that $P\setminus b$
contains no vertex of $B_1$.  This means that $a <_b a'$, which
contradicts the fact that $a\in\mbox{Max}(b)$.  So the lemma holds.
\hfill $\Box$

\medskip

Let $A_1=\{a_1, \ldots, a_k\}$ and $B_1=\{b_1, \ldots, b_\ell\}$, and
assume up to symmetry that $k\le\ell$.  By Lemma~\ref{lem:a1b1}, we
may assume up to relabeling that $\{a_1, b_1\}$ is an even pair of
$G$.  Similarly, for $i=2, \ldots, k$, we may assume that $\{a_i,
b_i\}$ is an even pair of $G\setminus\{a_1, b_1, \ldots, a_{i-1},
b_{i-1}\}$.

\subsection{Decomposing the graph}

By Lemma~\ref{lem:mabcuts}, the set $M\cup A_1\cup B_1$ is a cutset of
$G$, so $V(G)\setminus (M\cup A_1\cup B_1)$ can be partitioned into
two subsets $X$ and $Y$, with $C_1\subseteq X$ and $C_2\subset Y$,
such that there is no edge between $X$ and $Y$.  Let $G_X = G\setminus
Y$ and $G_Y = G\setminus X$.  Thus we consider that $G$ is decomposed
into $G_X$ and $G_Y$.  Since $G_X$ and $G_Y$ are proper induced
sugraphs of $G$, we may assume by induction that we have a clique
$Q_X$ of $G_X$ of size $\omega(G_X)$ and a coloring $c_X$ of $G_X$
with $\omega(G_X)$ colors, and the same for $G_Y$.
\begin{lemma}\label{lem:cx}
There exists a coloring $c'_X$ of $G_X$ with $\omega(G_X)$ colors such
that $c'_X(a_i)=c'_X(b_i)$ for all $i=1, \ldots, k$, and such a
coloring can be obtained from $c_X$ in polynomial time.
\end{lemma}
\noindent{\it Proof.} Suppose that $c_X$ itself does not have the
property described in the lemma, and let $h$ be the smallest integer
such that $c_X(a_h)\neq c_X(b_h)$.  In case $h>1$, we may assume, up
to relabeling, that $c_X(a_i)=i=c_X(b_i)$ for all $i=1,\ldots, h-1$.
Let $c_X(a_h)=i$ and $c_X(b_h)=j$, with $i\neq j$.  Note that both
$i,j> h-1$.  Let $H_{i,j}$ be the bipartite subgraph of $G$ induced by
the vertices of color $i$ and $j$.  We swap colors $i$ and $j$ in the
component of $H$ that contains $a_h$.  This component does not contain
$b_h$, for otherwise it contains a chordless odd path between $a_h$
and $b_h$, and this path is in $G\setminus\{a_1,b_1,\ldots,
a_{h-1},b_{h-1}\}$ since it contains no vertex of color less than $i$
and $j$; but this contradicts the fact that $\{a_h,b_h\}$ is an even
pair of $G\setminus\{a_1,b_1,\ldots, a_{h-1},b_{h-1}\}$.  So after
this swapping vertices $a_h$ and $b_h$ have the same color.  Thus we
obtain a coloring of $G_X$ with $\omega(G_X)$ colors where the value
of $h$ has increased.  Repeating this procedure at most $k$ times
leads to the desired coloring.  \hfill $\Box$

\medskip

Applying Lemma~\ref{lem:cx} to both $G_X$ and $G_Y$, we obtain
colorings $c_X$ and $c_Y$ of $G_X$ and $G_Y$ respectively such that,
up to relabeling, $c_X(a_i)=c_X(b_i)=c_Y(a_i)=c_Y(b_i)=i$ for each
$i=1, \ldots, k$.  Recall that $M\cup (B\setminus\{b_1, \ldots,
b_k\})$ is a clique and that all its vertices are adjacent to at least
one of $a_i, b_i$ for each $i=1,\ldots,k$.  So we may assume, up to
relabeling, that every vertex $z$ in $M\cup (B\setminus\{b_1, \ldots,
b_k\})$ satisfies $c_X(z)=c_Y(z)$ too.  It follows that the two
colorings $c_X$ and $c_Y$ can be merged into a coloring of $G$.  This
coloring uses $\max\{\omega(G_X),\omega(G_Y)\}$ colors, and one of
$Q_X$ and $Q_Y$ is a clique of that size.  So the coloring and the
larger of these two cliques are both optimal.

\subsection{The algorithm}

We can now describe our algorithm.  
\begin{quote}
Input: A  graph $G$ on $n$ vertices.  \\
Output: Either a coloring of $G$ and a clique of the same size, or the
answer ``$G$ is not a square-free Grenoble graph''.  \\
Procedure: \\
1.  First test whether $G$ is square-free, and test whether $G$ is
Berge with the algorithm from \cite{CCLSV}.  Then test whether $G$
contains a prism as explained in \cite{MTa}.  If these tests produce
an induced subgraph of $G$ that is either a square, or an odd hole, or
an odd prism, return the answer ``$G$ is not a square-free Grenoble
graph'' and stop.  If the algorithm from \cite{MTa} shows that $G$
contains no prism, then color $G$ applying the algorithm from
\cite{MT}.  \\
2.  Now suppose that $G$ contains an even prism.  Grow a maximal
hyperprism $\eta$, and find a good strip $S_1$ of $\eta$.  \\
Apply the proof of Lemma~\ref{lem:order} to every vertex $x\in A_1\cup
B_1$.  That proof either establishes that $<_x$ is an order relation
or finds $1$-rungs that converge; in the latter case, apply the proof
of Lemma~\ref{lem:conv} to obtain a new maximal hyperprism with a
smaller good strip, and restart from that hyperprism.   \\
When $<_x$ is an order relation for all $x\in A_1\cup B_1$,
Lemma~\ref{lem:a1b1} shows how to find even pairs.  The graph $G$ is
decomposed into graphs $G_X$ and $G_Y$, and an optimal coloring and a
maximal clique for $G$ can be obtained as explained above.
\end{quote}

Let us analyse the complexity of the algorithm.  One can decide
whether a given graph $G$ is Berge in time $O(n^9)$ with the algorithm
from \cite{CCLSV}.  One can test whether $G$ is square-free in time
$O(n^4)$, and whether a Berge graph $G$ contains a prism in time
$O(n^5)$ as explained in \cite{MTa}.  Now assume that the algorithm
produces an even prism.  It is easy to to see that all the procedures
in part 2 of algorithm (growing a maximal hyperprism, determining the
orderings) can be performed in time at most $O(n^3)$, and we make
additional remarks.  First remark that when we need to restart from a
new hyperprism, the size of the good strip is strictly smaller, and so
this restarting step occurs at most $O(n)$ times.  Secondly, remark
that when $G$ is decomposed into graphs $G_X$ and $G_Y$, the algorithm
is called recursively on them.  This defines a decomposition tree $T$
for $G$: every decomposition node of $T$ is an induced subgraph $G'$
of $G$ and has two children which are induced subgraphs of $G'$; and
every leaf of $T$ is a graph that contains no prism.  Let us show that
this tree has polynomial size.  When $G$ is decomposed into graphs
$G_X$ and $G_Y$ as above, because of a certain cutset that arises from
a hyperprism $\eta$, we mark the corresponding node of the tree with a
pair of vertices $\{c_1, c_2\}$ where $c_1\in C_1$ and $c_2\in C_2$
are chosen arbitrarily.  We mark every subsequent decomposition node
similarly.  Note that only pairs of non-adjacent vertices are used to
mark any node.
\begin{lemma}\label{lem:c1c2}
Every pair of vertices of $G$ is used to mark at most one node of the
decomposition tree.
\end{lemma}
\noindent{\it Proof.} Without loss of generality let us consider the
node $G$ itself, decomposed into graphs $G_X$ and $G_Y$ along a cutset
$M\cup A_1\cup B_1$ corresponding to a hyperprism $\eta$, with the
same notation as above.  Let $T_X$ be the subtree of $T$ whose root is
$G_X$, and define $T_Y$ similarly.  The node $G$ of $T$ is marked with
a pair of vertices $\{c_1, c_2\}$ where $c_1\in C_1$ and $c_2\in C_2$.
Since $c_1\notin Y$ and $c_2\notin X$, the pair $\{c_1,c_2\}$ is not
included in the vertex-set of any descendant of $G$ in the tree; so
this pair will not be used to mark any node of $T$ other than $G$.

Now suppose that a pair $\{c,d\}$ is used to mark a node in $T_X$ and
also a node in $T_Y$.  Then $\{c,d\} \subseteq V(G_X)\cap V(G_Y) =
M\cup A_1\cup B_1$, and since $c$ and $d$ are not adjacent, we have
$c\in A_1$ and $d\in B_1$.  Since $\{c,d\}$ marks a node in $T_X$,
there is a hyperprism $\eta_X$ in $G_X$ such that $c$ and $d$ lie in
the interior of two distinct strips of $\eta_X$.  Let $R_c$ and $R_d$
be rungs of $\eta_X$ that contain $c$ and $d$ respectively (so $R_c$
and $R_d$ lie in different strips of $\eta_X$), and let $R$ be a rung
in the third strip of $\eta_X$.  Let $K$ be the prism formed by $R_c,
R_d, R$.  So $V(K)\subseteq V(G_X)$.  Since $c\in A_1$, and $A_1$ is a
clique, and $c$ lies in the interior of $R_c$, it follows that $A_1$
contains at most one corner of $K$.  Likewise $B_1$ contains at most
one corner of $K$.  This implies that the set of major neighbors of
$K$ is included in $G_X$.  Moreover, if $A_1$ contains a corner $u$ of
$K$ and $B_1$ contains a corner $v$ of $K$, then $u$ and $v$ are not
in the same rung of $K$ (for otherwise $c$ and $d$ would also lie on
that same rung).  Let $R_2$ be any $2$-rung in $\eta$.  Then $R_2$
contains no major neighbor of $K$, and $R_2$ satisfies the hypothesis
of Theorem~\ref{spgt101} with respect to $K$.  The preceding
observations imply that $R_2$ must satisfy item~1 of
Theorem~\ref{spgt101}, and consequently $G$ contains an odd prism, a
contradiction.  So one of $T_X, T_Y$ is such that none of its nodes is
marked with $\{c,d\}$.  (Actually the preceding argument holds for
$T_Y$ as well, so any pair $\{c,d\}$ with $c\in A_1$ and $d\in B_1$
will never be used to mark any node of $T$.)

The preceding two paragraphs, repeated for every node of $T$, imply
the validity of the lemma.  \hfill $\Box$

\medskip

By Lemma~\ref{lem:c1c2} the total number of nodes in $T$ is $O(n^2)$.
The leaves of the decomposition tree $T$ are Berge graphs with no
antihole (since they are square-free) and no prism, so they can be
colored in time $O(n^6)$ as explained in \cite{MT}.  At each node $G'$
of $T$ different from the root $G$, we know that $G'$ is an induced
subgraph of $G$, so it is square-free Berge; hence we must only test
whether $G'$ contains a prism, which is done in time $O(n^5)$ as
explained in \cite{MTa}.  So the total complexity of the algorithm is
$O(n^9+ n^2\times n^5 + n^2\times n^6)$ $=O(n^9)$.  This completes the
proof of Theorem~\ref{thm:main}.

\section*{Acknowledgement}

The author is partially supported by ANR project STINT under reference
ANR-13-BS02-0007.


\end{document}